\documentclass[leqno, 12pt]{amsart}
\usepackage{amssymb,amsmath,amsthm,graphicx,multirow,setspace,url}
\usepackage{amscd}
\usepackage{tikz-cd}
\usepackage{mathrsfs}

\usepackage{graphicx}

\usetikzlibrary{arrows}
\usetikzlibrary{quotes}

\usepackage{lipsum}% http://ctan.org/pkg/lipsum
\usepackage{fancyhdr}% http://ctan.org/pkg/fancyhdr
\fancypagestyle{mypagestyle}{%
  \fancyhf{}% Clear header/footer
  \fancyhead[OC]{Chow groups, Riemann-Roch without Denominators}% Author on Odd page, Centred
  \fancyhead[EC]{S. P. Dutta}% Title on Even page, Centred
  \fancyhead[L]{\thepage}
}
\title{A title}
\title{Chow groups: A Structure Theorem, Riemann-Roch without Denominators and Artin Approximation}
\author{S. P. Dutta}
\address{Department of Mathematics\\
University of Illinois at Urbana-Champaign\\
1409 West Green Street\\
Urbana, Illinois 61801}

\begin{document}

\newcommand{\hookuparrow}{\mathrel{\rotatebox[origin=c]{90}{$\hookrightarrow$}}}
\newcommand{\hookdownarrow}{\mathrel{\rotatebox[origin=c]{-90}{$\hookrightarrow$}}}

\newcommand{\simequp}{\mathrel{\rotatebox[origin=c]{90}{$\simeq$}}}
\newcommand{\simeqdown}{\mathrel{\rotatebox[origin=c]{-90}{$\simeq$}}}

%\numberwithin{equation}{section}
%\numberwithin{table}{section}

%\theoremstyle{custom}
\theoremstyle{plain}
\newtheorem*{theorem*}{Theorem}
\newtheorem{theorem}{Theorem}

\newtheorem{innercustomthm}{Theorem}
\newenvironment{customthm}[1]
  {\renewcommand\theinnercustomthm{#1}\innercustomthm}
  {\endinnercustomthm}

%newtheorem*{theorem}{Theorem}
%
%\newtheoremstyle{named}{}{}{\itshape}{}{\bfseries}{.}{.5em}{\thmnote{#3's }#1}
%\theoremstyle{named}
%\newtheorem*{namedtheorem}{Theorem}

%\newtheorem{theorem}{Theorem}[section]
%\newtheorem*{theorem*}{Theorem}[section]
%\newtheorem{theorem}{Theorem}

%\begin{theorem}  A numbered theorem.    \end{theorem}
%\begin{theorem*} An unnumbered theorem. \end{theorem*}

\newtheorem*{corollary}{Corollary}

\theoremstyle{remark}
\newtheorem*{remark}{Remark}
\newtheorem*{claim}{Claim}

\newtheorem*{namedtheorem}{\theoremname}
\newcommand{\theoremname}{testing}
\newenvironment{named}[1]{\renewcommand{\theoremname}{#1}
  \begin{namedtheorem}}
	{\end{namedtheorem}}
%New function names and their shortcuts

\newcommand{\gr}{\operatorname{grade}}
\newcommand{\Syz}{\operatorname{Syz}}
\newcommand{\syz}[2]{\Syz^{#1}(#2)}
\newcommand{\ring}[2]{{\mathcal{O}_{#1}(#2)}}
\newcommand{\hm}[3]{H_{#1}^{#2}(#3)}
\newcommand{\Tor}{\operatorname{Tor}}
\newcommand{\tor}[3]{\Tor_{#1}^{#2}(#3)}
\newcommand{\Ext}{\operatorname{Ext}}
\newcommand{\ext}[3]{\Ext_{#1}^{#2}(#3)}
\newcommand{\Id}{\operatorname{Id}}
\newcommand{\im}{\operatorname{Im}}
\newcommand{\coker}{\operatorname{coker}}
\newcommand{\grade}{\operatorname{grade}}
\newcommand{\Ht}{\operatorname{height}}
\newcommand{\Hom}{\operatorname{Hom}}
\newcommand{\Der}{\operatorname{Der}}
\newcommand{\dm}{\operatorname{dim}}

%shortcuts
%\newcommand{\mc}[1]{\mathcal{#1}}
\newcommand{\ul}[1]{\underline{#1}}

\begin{abstract}
The focus of this note is on the Chow group problem over ramified regular local rings $(R, m)$. Our goal is threefold: i) to introduce a characterization of a ramified regular local ring essentially of finite type over a dvr, ii) to address the question whether $(i-1)!$ $\mathbb{A}^i(U)=0$ for specific open subsets $U$ of Spec$R$ and iii) to establish a constructive relation between Chow groups of the henselization $(R^h, m^h)$ and Chow groups of the completion $(\hat{R}, \hat{m})$ of $(R, m)$.
\end{abstract}

\maketitle

\section{Introduction}\label{Intro}
\textbf{1.1.} The Chow group problem for any regular local ring $(R, m)$ of dimension $n$ asserts that the Chow groups $\mathbb{A}_i(R)=0$ for $0\leq i < n$. Let $\mathbb{A}^i(R)=\mathbb{A}_{n-i}(R)$. Claborn and Fossum raised the following question in [C-F]: Given any regular local ring $(R, m)$ of dimension $n$, is $(i-1)!$ $ \mathbb{A}^i(R)=0$? Using generalized Riemann-Roch theorem (20.1, [Fu]) Fulton proved that over any regular local ring Chow groups are torsion-groups (20.1.4, [Fu]); no precise integer, as above, for torsion was suggested (see brief history). In this note our main goal is threefold: 1) to introduce a characterization of a ramified regular local ring essentially of finite type over a dvr, 2) to address this question on specific open subsets of spec$R$ for a ramified regular local ring $R$  and 3) to establish constructive relation between $\mathbb{A}^i({R}^h)$ and $\mathbb{A}^i(\hat{R})$ and their corresponding $K_0$ groups, where $R^h$ and $\hat{R}$ are henselization and completion of a ramified regular local ring $R$ respectively with respect to the maximal ideal $m$. Let us mention the main results briefly. 

\smallskip
In section 2 we prove the following characterization of a ramified
\smallskip

\line(1,0){300}

AMS Subject Classification:

Primary:  13D22, 14C40;
Secondary: 13H05

Key words and phrase: Regular local ring, Chow group, Riemann-Roch theorem without denominators, Artin approximation.

\noindent
regular local ring $(R, m)$ essentially of finite type over a pseudo-geometric discrete valuation ring, henceforth dvr, $(V, pV)$ in mixed characteristic $p>0$. Recall that a regular local ring $(R, m)$ in mixed characteristic $p$ ($p=ch.R/m$) is called ramified if $p \in m^2$ and a Noetherian ring $A$ is called pseudo-geometric if for any prime ideal $P$ of $A$, the integral closure of $A/P$ in any finite extension of the quotient field of $A/P$ is a finitely generated $A/P$ module.

\smallskip
\textbf{Theorem(2.1.)} Let $(R, m)$ be as above and let dimension of $R$ be $n$. We have the following:

\smallskip
i) Assume that $V/pV$ is perfect. Then $R \simeq (S[X]/(f(X)))_{\tilde{m}}$, where $S = W[X_1, ..., X_{n-1}]_{(p, X_1, ..., X_{n-1})}$, $(W, pW)$ is a pseudo-geometric dvr, $f(X)$ is monic irreducible in $S[X]$ and $\tilde{m}$ is a maximal ideal of \\ $S[X]/(f(X))$.

ii) Suppose $V/pV = R/m$. Then $R \simeq (S[X]/(f(X)))_{\tilde{m}}; S, f(X), \tilde{m}$ as above and $X \in \tilde{m}$.

Note that case ii) of the above theorem resembles the structure theorem for complete  ramified regular local rings due to Nagata (Th.31.12, [N]). Theorem 2.1 would play an important role in our study of the Chow group problem in this note.
\smallskip

In section 3 first we deal with smooth varieties in the equicharacteristic case. We have the following.
\smallskip

\textbf{Theorem(3.1)} Let $X$ be a smooth variety of dimension $n$ over a field $k$. Assume that $K_0(X) = \mathbb{Z}$, $K_0(X) = K_0(\mathscr{M}_0)$ where $\mathscr{M}_0$ is the category of coherent sheaves on $X$. Then $(i-1)!$ $\mathbb{A}^i (X) = 0$ for $0<i\leq n$.

\smallskip
Riemann-Roch without denominators in equicharacteristic over a field $k$ (eqn.(15), p.150, eqn.(16), p.151, [Gr1]; cor.2, p.64, [Gr2]; sec.15.3 [Fu]) plays a crucial role in our proof. Grothendieck introduced the concept of Riemann-Roch without denominators and mentioned the result in our paper (cor.4.2) in (eqn.(16), p.151, [Gr1]). Jouanalou [J] gave a proof of Grothendieck's conjecture on Riemann-Roch without denominators for closed imbeddings of smooth varieties (expos\'e XIV, GR2). Fulton sketched a different proof for the smooth case (15.3, [Fu]).
\smallskip

We derive the following corollaries.
\smallskip

\textbf{Corollary (3.2).} Let $R$ be a ramified regular local ring of dimension $n+1$ essentially of finite type over a dvr $(V, pV)$. Then $(i-1)!$ $\mathbb{A}^i (R[1/p])=0$ for $0<i \leq n$.
\smallskip

\textbf{Corollary(3.3).} Let $(R, m)$ be as in Cor.(3.2) and let $R^h$ denote its hensalization. Then $(i-1)! \mathbb{A}^i(R^h[1/p])=0$ for $0<i \leq n$. 
\smallskip

\textbf{Corollary(3.4)} Let $(R, m)$ be as above. Let $\bar R = R/pR$; for any prime ideal $P$ of $R$, $p \notin P$, let $\bar R/\bar P= R/(P+pR)$. If codimension of $P$ is $i$, then $(i-1)![\bar R/\bar P]=0$ in $\mathbb{A}_{n-i}(\bar R)$.

\smallskip

In section 4 we want to prove similar results as in section 3 for another important open subset of spec$R$ in mixed characteristic $p>0$. Due to theorem 2.1, if $R$ is a ramified regular local ring in mixed characteristic $p>0$ essentially of finite type over a dvr ($V, pV$) such that $V/pV$ is perfect, then $R=S[X]/(f(X)) \tilde{m}$, where $S, f(X)$ and $\tilde{m}$ are as stated in theorem 2.1. Hence $V \rightarrow R_{f'}$ is smooth. In order to prove results parallel to cor.(3.2), cor.(3.3) and cor.(3.4), via our approach, we need to prove an extended version of theorem 3.1 in mixed characteristic. This led us to look for a mixed characteristic version of Riemann-Roch theorem without denominators. After checking sketches of proofs in [(GR1), (GR2) and (Fu)] we realised that Fulton's arguments in [(Fu)] can be extended to prove the required version in mixed characteristic. Unable to find any reference in the literature about such an extension, we asked Fulton whether Riemann-Roch theorem without denominators should be valid in mixed characteristic. His reply confirmed our observation. 

\smallskip

In section 4.1 for the benefit of the readers we briefly mention the apparent road blocks for such an extension in any characteristic for smooth varieties over a Dedekind domain and point out Fulton's work in Chapter 20, [(Fu)] (preceded by works of Grothendieck, Samuel and Shimura) to overcome these obstacles. In 4.2 we propose the following theorem.

\smallskip

%In section 4 first we want to present a mixed characteristic analog of Riemann-Roch theorem without denominators in equicharacteristic. This observation seems to be new as it has not been mentioned in chapter 20 [Fu]. When we first realized that this mixed characteristic case should be valid (e.g., replacing field by dvr), we asked Fulton whether Riemann-Roch theorem without denominators under right conditions (see below) is valid in mixed characteristic. A part of his reply was the following "I certainly think so and expect that there is no problem in verifying that the proof for a field extends to this case." This is exactly true. In this section we mention the apparent road blocks in mixed characteristic and point out Fulton's work in chapter 20, [Fu] to overcome this obstacles. We propose the following theorem and a corollary.

\smallskip

\textbf{Theorem(4.2).} (Riemann-Roch theorem without denominators in any characteristic [expos\'e XIV, GR2], [J], 15.3.1 [Fu]).

\smallskip

Let $S=$Spec$A$ where $A$ is a Dedekind domain in any characteristic. Let $f:X \rightarrow Y$ be a closed imbedding of smooth varieties over Spec $S$ of codimension $d$ with normal bundle $N$. Let $E$ be a vector bundle of rank \ $e$ on $X$. Then 
$$c(f_*[E])=1+f_*(g(N, E))$$
where $g(N, E)$ is a unique power series with integer coefficients in variables $T_1, ..., T_d, U_1, ..., U_e$ (for description of $g$ we refer the reader to \textbf{observation 1} in section 4.2).

\smallskip

We should mention here that for our sketches of proofs of {\bf observation 1} and the above theorem in section 4.2, we offer arguments for extensions of crucial steps in Fulton's work (Lemma 15.3, theorem 15.3.1, [Fu]) leaving the details to the reader. 

\smallskip

We have the following important corollary.

\smallskip 

\textbf{Corollary(4.3).} With assumption as above we have the following:

$c_j(f_*[O_X])=0, 0<j<d, c_d (f_* [O_X])= (-1)^{d-1} (d-1)! [X]$ in $\mathbb{A}^d (Y/S)$ (GR1, 15.3.2 [Fu], [P-Sh]).

A proof of this corollary follows directly from the above theorem. Actually we provide a direct proof of this corollary for all characteristics without appealing to the Riemann-Roch theorem without denominators mentioned in theorem 4.2.

Our final theorem in this section is a similar version of theorem 3.1 over a Dedekind domain in any characteristic.

\smallskip
\textbf{Theorem(4.4).} Let $X$ be a smooth variety of dimension $n+1$ over a Dedekind domain $A$ in any characteristic. Assume that $K_0(X)=\mathbb{Z}$, where $K_0(X)=K_0(\mathscr{M}_0), \mathscr{M}_0$ is the category of coherent sheaves on $X$. Then $(i-1)! \mathbb{A}^i (X)=0$ for $0<i \leq n+1$.
\smallskip

We have the following corollaries.
\smallskip

\textbf{Corollary(4.5).} Let $(R, m)$ be a ramified regular local ring of dimension $(n+1)$ essentially finite type over a pseudo-geometric dvr $(V, pV)$ such that $V/pV$ is perfect or $V/pV = R/m$. We have $R=(S[X]/(f(X)))_{\tilde{m}}$ by theorem 2.1. Then the unramified locus of Spec$R$ over Spec$S$, i.e., Spec$R_{f'}$ has the property that, for $0<i \leq n, (i-1)!$ $\mathbb{A}^i(R_{f'})=0$.
\smallskip

\textbf{Corollary(4.6).} Let $(R, m)$ be as above and let $(R^h, m^h)$ be its henselization. Let $U^h=$ quasi-unramified locus of Spec$R^h$ over Spec$S$. Then $(i-1)!$ $\mathbb{A}^i(U^h) =0$ for $0<i \leq n$.
\smallskip

\textbf{Corollary(4.7).} $(R, m)$ is as in the corollary above. Let $P$ be a prime ideal of $R$ such that ht$\bar{P=i}$ and $f' \notin P$. Let $\bar R=R/f'R$ and $\bar R/\bar P = R/(P+f' R)$. Then $(i-i)!  [\bar R/\bar P]=0$ in $\mathbb{A}^i (\bar{R})$.
\smallskip

In section 5 we deal with the relation between $\mathbb{A}^i (R^h)$ and $\mathbb{A}^i(\hat{R})$, where $(R^h, m^h)$ and $(\hat{R}, \hat{m})$ are the henselization and completion of a ramified regular local ring $(R, m)$ of dimension $n$ essentially of finite type over an excellent dvr respectively with respect to the maximal ideal $m$. Let us recall that in general the Chow group map induced by a flat morphism between two Noetherian schemes of finite type over a field or a dvr is neither injective nor surjective. Moreover, although $\hat{R}$ has a nice structure ([N]), $R^h$ does not, even when $R$ has a nice structure (Th. 2.1). Our effort to understand Chow groups of $R^h$ via direct limit process has been described in cor.(3.3) and cor.(4.6). Now we would like to investigate this problem via the flat map $R^h \hookrightarrow \hat{R}$. Over the years we tried to understand the effect of the validity of the Chow group problem over $\hat{R}$ on the validity of the same over $R^h$. Finally we are able to figure out the relation and it is presented in theorem 5.4 below. Lemma 5.2 and proposition 5.3 pave the way for our proof of this theorem.

\smallskip

In this section we utilize the presentation of Chow groups by Claborn \& Fossum in ([C-F]). Note that the definition of $W_i(R)$ as prescribed in (p.230, [C-F]) converges with the definition of $\mathbb{A}^i(R)$ over any regular local ring (actually over any Cohen-Macaulay local ring). Artin approximation plays an important role here. Results derived in [D3] using Artin approximation ([A]) and result proved in [D1] are used in our proofs of the results below (sec.5 for more details). 

\smallskip

Recall that for any regular local ring $A$, $M_i(A)$ denotes the Serre-subcategory of finitely generated $A$-modules of codimension$\geq i$. First we prove the following.

\smallskip

{\bf Lemma (5.2).} $\mathbb{A}^i(R)=0$ for $1 \leq i \leq n \iff K_0(\mathscr{M}_i (R)) \xrightarrow{\phi_{{i-1}, i}} K_0(\mathscr{M}_{i-1} (R))$ are the {0}-maps for all $i, 1 \leq i \leq n$. This is valid for any regular local ring.

Note that this shows that the validity of $K_0$-case of Gersten's conjecture [Q] on any regular local ring is equivalent to the validity of the Chow group problem.

\smallskip

{\bf Proposition (5.3).} $\forall i, 1\leq i \leq n, K_0(\mathscr{M}_i (R^h)) \rightarrow K_0(\mathscr{M}_i(\hat{R}))$ is injective.

\smallskip

Results 1 and 2 mentioned in section 5.1 pave the way for proving this proposition. 

\smallskip

Now we state our theorem and a corollary.

\smallskip

{\bf Theorem (5.4).} If $\mathbb{A}^j (\hat{R}) =0$ for $i \leq j \leq n$, then $\mathbb{A}^j (R^h)=0$ for $i \leq j \leq n$.

\smallskip

{\bf Corollary (5.5).} $\mathbb{A}^{n-1}(R^h)=0$.

\smallskip

%Let $\mathscr{M}_i$ denote the Serre-subcategory of finitely generated $R$-modules of codimension $\geq i$.

%We have the following theorem:

%\smallskip

%\textbf{Theorem(5.1).} With notations as above we prove the following:

%i) $\mathbb{A}^i(R)=0$ for $1 \leq i \leq n \iff K_0(\mathscr{M}_i) \xrightarrow{\phi_{{i-1}, i}} K_0(\mathscr{M}_{i-1})$ are the {0}-maps for all $i, 1 \leq i \leq n$. This is valid for any regular local ring.

%ii) $\forall i, 1\leq i \leq n, K_0(\mathscr{M}_i (R^h)) \rightarrow K_0(\mathscr{M}_i(\hat{R}))$ is injective.

%iii) If $\mathbb{A}^j (\hat{R}) =0$ for $i \leq j \leq n$, then $\mathbb{A}^j (R^h)=0$ for $i \leq j \leq n$.

%iv) $\mathbb{A}^{n-1}(R^h)=0$.
\medskip

\textbf{1.2 Brief History.} Fossum and Claborn proved in [C-F] that when $R$ is a power series ring/polynomial ring over a field $k$ or a discrete valuation ring $V, \mathbb{A}_i(R)=0$ for $i<n, n$=Krull dimension of $R$. $(\mathbb{A}_n(R)= \mathbb{Z}$ and $\mathbb{A}_0(R)=0$ are obvious). Quillen [Q] proved Gersten's conjecture for regular semi-local rings smooth over a field $k$. Fulton [Fu] pointed out that for any regular local ring the Chow groups are torsion groups. Later Gillet and Levine [G-L] used Quillen's techniques to extend the validity of Gersten's conjecture [Q] to the case when the above rings are smooth over a dvr. Using K-theoretic techniques, Levine [Lev] also proved the same result over $R=V[[X_1, ..., X_n]]/(p+X_1^2+...+X_n^2), p \neq 2$, $(p)=m_V$ and $p=$ mixed characteristic of $R$. The $K_0$-case of Gersten's conjecture [Q] is the Chow group problem for regular local rings (Lemma 5.2). The following results were proved in [D1] and [D2] for complete ramified regular local rings $(R, m)$ of dimension $n$. 

1. $\mathbb{A}_1(R)=0$. (Theorem 1.1, [D1])

2. Let $R=S[X]/(f(X))$, where $S$ is a power series ring in $(n-1)$ variables over a dvr $(V, p)$, $f(X)$ is an Eisenstein polynomial of degree $t$, $p$= mixed characteristic of $R$. Let $P$ be a prime ideal of $R$ of height $i$ and let $q=S \cap P$. 

a) If $S/q$ is normal then, by inducting on height of $P$, $[R/P]=0$ in $\mathbb{A}^i(R)$.

b) Suppose $R/P$ is normal. Two cases: i) $R/P$ is the integral closure of $S/q$ and ii) $R/P$ is not so. If degree of $f(X)$ is prime, then $[R/P]=0$ in $\mathbb{A}^i(R)$ in case (ii); when degree of $f(X)$ is not a prime, then validity of case (i) for the Chow group problem implies the validity of case (ii) for the same. Moreover a new proof of the equicharacteristic/unramified (restricted) case was proved in [D2] via an approach completely different from that due to Quillen. An alternative proof of Levine's theorem for Chow groups was also pointed out in [D2]. In his thesis at University of Illinois S. Lee [Le] solved the Chow group problem for complete regular local rings of dimension $\leq 4$. It is clear from above that the ramified case of the Chow-group problem is still very much open even in the complete case.

\medskip

\begin{center}{Section 2: A characterization of ramified regular local rings}
\end{center}
\smallskip

We prove the following theorem characterizing ramified regular local rings essentially of finite type over a dvr. This will be useful in our study of Chow groups.

\medskip

\textbf{2.1 Theorem.} Let $(V, p)$ be a pseudo-geometric dvr in mixed characteristic $p>0$. Let $(R, m)$ be a ramified regular local ring of dimension $n$ essentially of finite type over $V$. We have the following.
\smallskip

i) Assume that $V/pV$ is perfect. Then there exists a pseudo-geometric dvr $(W, p)$ such that if $S=W[X_1, ..., X_{n-1}]_{(p, X_1, ..., X_n)}$, then \\ $(S[X]/(f(X)))_{\tilde{m}} \simeq R$ where $\tilde{m}$ is a maximal ideal in $S[X]$ and $f(X)$ is monic irreducible in $S[X]$.
\smallskip

ii) Assume $V/pV \simeq R/m$. Then $R \simeq (S[X]/(f(X)))_{\tilde{m}}: S, f(X), \tilde{m}$ as above and $X \in \tilde{m}$. 
\medskip

\textbf{Proof.} i) Let $k=V/pV$. Since $k$ is perfect, $R/m$ is separately generated over $k$. Let $t_1, ... t_d$ be a separating transcendence basis of $R/m$ over $k$, i.e., $t_1, t_2, ...,., t_d$ are algebraically independent over $k$ and $R/m$ is a simple extension over $k(t_1, t_2, ..., t_d)$, i.e., $R/m=k(t_1, t_2, ..., t_d)({\alpha})$, $\alpha$ is separable algebraic over $k(t_1, t_2, ..., t_d)$. Let $T_1, T_2, ..., T_d$ denote lifts of $t_1, ..., t_d$ in $R$. Then $T_1, ..., T_d$ are algebraically independent over $V$. Let $W=V[T_1, ..., T_d]_{pV[T_i, ..., T_d]}$; then $W$ is a pseudo-geometric dvr with maximal ideal generated by $p$ such that $W/pW=k(t_1..., t_d)$. Let $X_1, ..., X_{n-1}$ be a part of a regular system of parameters of $R$ such that $p, X_1, ..., X_{n-1}$ form a system of parameters of $R$. Let $A=W[X_1, ..., X_{n-1}]$; then ht$(p, X_1, ..., X_i)=i+1$, $0 \leq i < n$ and $X_1, ..., X_{n-1}$ are algebraically independent over $W$. By the dimension formula, $R$ is algebraic over $A$. Let $B$ denote the integral closure of $A$ in $Q(R)$. Since $R$ is regular local, $B \subset R$ and $B$ is a module finite extension of $A$. Let $S=W[X_1, ..., X_{n-1}]_{(p, X_1, ..., X_{n-1})}$ and $T=B_{(p, X_1, ..., X_{n-1})}$. Then $T$ is a semi-local module-finite extension of $S$. Let $j:T \hookrightarrow R$ denote the natural injection. Let $m_1, ..., m_r$ denote the maximal ideals of $T$ where $m_1=j^{-1}(m)$. Then, due to an observation by Nagata (Th. 37.4, [N]), $T_{m_1} = R$, $T/m_1=R/m = S/m_S ({\alpha})$ where $m_S$ denotes the maximal ideal of $S$, $m_S=(p, X_1, ..., X_{n-1})S$. We have
\smallskip

$T/m_S T \simeq T/I_1 \times ... \times T/I_r$ ...................(1) 
\smallskip

where $I_i$ is $m_i$-primary ideal in $T$ for $1 \leq i \leq r$. Since $R$ is ramified, i.e., $p$ is in $m^2$, $I_1=(p, X_1, X_2, ..., X_{n-1})R \varsubsetneq m$; moreover, since $p \in m^2$, the maximal ideal $m_S$ of $S$ generates the subspace spanned by im$X_1, ...,$ im$X_{n-1} \in m/m^2=m_1/m_1^2$. Let $\beta \in T$ be such that the image of $\beta = \alpha \in T/m_1$. Let $g(X)$ denote a monic polynomial in $S[X]$ whose reduction mod $m_S$ is the minimal polynomial for $\alpha$. Then
\medskip

$g(\beta) \in m_1, g'(\beta) \notin m_1$ and $g(0)$ is a unit in $S ..........(2).$ 
\medskip

We can choose this $\beta$ in such a way that $X_1, X_2, ..., X_{n-1}, g(\beta)$ form a regular system of parameters of $R$, i.e., a minimal set of generators of $m$; equivalently, a lift of a basis of $m/m^2 = m_1/m_1^2$. This can be done, if needed, by replacing $\beta$ with $\beta + h$, where $X_1, ..., X_{n-1}$, $h$ form a regular system of parameters of $R$ (by (2)), i.e., $\{\bar{X_1}, ..., \bar{X_{n-1}}, \bar{h} \}$ form a basis of $m_1/m_1^2= m/m^2$ over $T/m_1=R/m$. Thus $X_1, ..., X_{n-1}, g(\beta)$ in $T$ form a regular system of parameters for $T_{m_1}=R$. 

\smallskip

Let $\gamma \in T$ be such that $\gamma$ lifts $(\bar{\beta}, 0, ..., 0)$ in $T/m_S T$ where $\bar{\beta}=$ im $\beta$ in $T/I_1$ (1). Then $g(\gamma)$ does not belong to $m_i$ for $2 \leq i \leq r$ (by (2)) and $X_1, ..., X_{n-1}, g(\gamma)$ form a regular system of parameters for $T_{m_1}=R$. We now consider the ring $C=S[\gamma]$. Let $m'_i$ denote the maximal ideal $m_i \cap C$ for $1 \leq i \leq r$. We write $q=m'_1$. Since $g(\gamma) \notin m_i, 2 \leq i \leq r$, no $m_i$ in $T$ can lie over $q$ in $C$ for $2 \leq i \leq r$. Thus $m_1$ is the only maximal ideal in $T$ lying over $q$ in $C$. Hence $T_q = T_{m_1}=R$. Since $T_q$ is a finitely generated $C_q$ module and $C/q$ is isomorphic to $T/m_1$ by construction, we have $C_q=T_{m_1}=R$.

\smallskip

$\gamma$ is integral over $S$; let $f(X)$ denote the minimal polynomial of $\gamma$ in $\mathbb{Q}(S[X])$. Since $S$ is normal, $f(X) \in S[X]$. Hence $S[X]/(f(X)) \simeq C$, via the map $X \rightarrow \gamma$. Let $\tilde{m}$ be the inverse image of $q$ in $S[X]$. Then $R=C_q=(S[X]/(f(X)))_{\tilde{m}}$.

\smallskip

\textbf{ii)} If $V/pV \simeq R/m$, then $g(\gamma)=\gamma-a, a \notin m_S$. This implies $X-a \in \tilde{m}$. By means of a change of variable we obtain $R=(S[X]/(f(X)))_{\tilde{m}}$, where $X \in \tilde{m}$. 
\smallskip

\textbf{Remark.} Part ii) of the above theorem resembles Nagata's characterization of a complete ramified regular local ring (Th. 31.12, [N]). However, we do not know whether $f(X)$ is an Eisenstein polynomial. We can only assert that the constant term of $f(X)$ and coefficient of $X$ in $f(X)$ must be in $m_S$. 
\medskip

\begin{center}{Section 3: Equicharacteristic case over a field and its application.}
\end{center}
\smallskip

We have the following theorem.

\smallskip

\textbf{3.1 Theorem.} Let $X$ be a smooth variety of dimension $n$ over a field $k$. Assume that $K_0(X)=\mathbb{Z}$; $K_0(X)$ denotes $K_0(\mathscr{M}_0)$, where $\mathscr{M}_0$ is the category of coherent sheaves on $X$. Then $(i-1)!$ $\mathbb{A}^i(X)=0$ for $0 <i \leq n$.
\medskip

\textbf{Proof.} Let $\mathscr{M}_i$ be the Serre subcategory of coherent sheaves of codimension $\geq i$. We have a descending chain
$$\mathscr{M}_0 \supset \mathscr{M}_1 \supset ... \supset \mathscr{M}_i \supset \mathscr{M}_{i+1} \supset ...$$
and a corresponding homomorphism of $K_0$-groups for $i \leq j, \phi_{ij}:K _0 \\(\mathscr{M}_j) \rightarrow K_0(\mathscr{M}_i)$. Let $F^i=\phi_{0i}(K_0(\mathscr{M}_i)); \{F^i\}_{i \geq 0}$ constitute a filtration on $K_0(X)$. Let $G(K_0(X))$ denote the associated graded group, i.e., $G_i(K_0(X))=F^i/F^{i+1}$. Recall that $\mathbb{A}^i(X)=\mathbb{A}_{n-i}(X)$. It has been pointed out in ([C-F], [Fu], [Gr1], [Gr2]) that, for each $i$, there exists group homomorphism

$$\phi_i:\mathbb{A}^i(X) \rightarrow F^i/F^{i+1}=G_i(K_0(X))$$ given by $\phi_i([V]) = $ class of $O_V$ for any sub-variety $V$ of codimension $i$ and extending by linearity it induces a graded group homomorphism 
\smallskip

$\phi: \mathbb{A}^*(X) \rightarrow G(K_0(X))$ where $\mathbb{A}^*(X)=\mathbb{\oplus} \mathbb{A}^i(X)$. Recall that since $X$ is smooth over $k, \mathbb{A}^*(X)$ also has a multiplicative structure:
\smallskip

$\mathbb{A}^i(X) \otimes \mathbb{A}^j(X) \rightarrow \mathbb{A}^{i+j}(X)$ given by $x \otimes y \rightarrow x.y$, the corresponding intersection product. $\phi$ is surjective (and commutes with proper push-forward).

On the other hand while introducing Riemann-Roch theorem without denominators on smooth quasi-projective varieties over $k$, Grothen-dieck defined (later Fulton defined on smooth varieties over $k$) a graded homomorphism ([Fu], [Gr1], [Gr2]) 

$\psi:GK_0(X) \rightarrow \mathbb{A}^*(X)$ in the following way. Let $E$ be a vector bundle on $X$; let $c(E)= 1+c_1(E)+ ... $ denote the total Chern class as an element of the multiplicative group $\mathbb{A}^*(X)$. If $E_{\bullet}$ is a complex of vector bundles on $X$, $c(E_{\bullet})=c(\Sigma(-1)^i[E_i]) = \Pi c(E_i)^{(-1)^i}$. If $\alpha \in F^i(X)$, then $c_j(\alpha)=0$ for $0<j<i$ and hence $\psi_i:F^i/F^{i+1} \rightarrow \mathbb{A}^i(X)$ is defined by $\psi_i$ (class of $\alpha) = c_i(\alpha) \in \mathbb{A}^i(X)$. Thus we get a group homomorphism $\psi:GK_0(X) \rightarrow \mathbb{A}^*(X)$.

It has been mentioned in ([Gr1]) and proved in ([Fu], [Gr2]) that for each $i$, the composite maps $\phi_i \bullet \psi_i$ and $\psi_i \bullet \phi_i$ satisfy 

\smallskip
$\phi_i \bullet \psi_i = (-1)^{i-1}(i-1)! \ id$ and $\psi_i \bullet \phi_i = (-1)^{i-1} (i-1)! \ id.................(3)$
\smallskip

Let us explain the above facts a bit more elaborately in the following steps.

a) $Y$ - a closed subscheme of $X$. Due to the exact sequence 
$$\mathbb{A}_t(Y) \rightarrow \mathbb{A}_t(X) \rightarrow \mathbb{A}_t (X-Y) \rightarrow 0, ........(4)$$
It follows that if dim$Y<t$, then $\mathbb{A}_t(X) \simeq \mathbb{A}_t (X-Y)$.

b) Let $\alpha = O_W$, codimension of $W=i$, i.e., $\alpha \in F^i (X)$. Let $f:W \hookrightarrow X$ denote the natural injection. Consider a resolution of $(f_* (O_W))$ by locally free sheaves $E_{\bullet}$ over $X$. Then 
$$c(f_* [O_W]) = c(\Sigma (-1)^i [E_i]) = \Pi c(E_i)^{(-1)^i}.$$
We have for any $j$ positive, $c_j (f_* [O_W]) =0$ on $X-W$. Hence, replacing $Y$ by $W$, it follows from (4) that for $j<i$, $c_j (f_* [O_W])=0$ and $c_i (f_* [O_W])=$ a multiple of [$W$] in $\mathbb{A}^i (X)$. The main point is to show 
\smallskip

$c_i (f_* [O_W]) = (-1)^{i-1} (i-1)! [W]$ in $\mathbb{A}^i (X) ......(5)$.
\smallskip

When $W$ is smooth, (5) follows from the main crux of the Riemann-Roch theorem without denominators and the reader is referred to theorem 15.3 and example 15.3.1 in [Fu]. When $W$ is not smooth, let $Z$ denote the non-smooth locus of $W$. Then $W-Z$ is smooth and $W-Z \hookrightarrow X-Z$ is a regular imbedding. Since codim $Z>i, \mathbb{A}^i (X) \simeq \mathbb{A}^i (X-Z)$ by (4). Hence, in this case (5) follows from the above arguments.

\smallskip

By assumption, we have $K_0(X)=\mathbb{Z}$. We have a short exact sequence
$$K_0(\mathscr{M}_1) \xrightarrow{\phi_{01}} K_0(\mathscr{M}_0) \rightarrow K_0(\mathscr{M}_0/\mathscr{M}_1) \rightarrow 0$$

Since $K_0(\mathscr{M}_0/\mathscr{M}_1)=\mathbb{Z}$, we have $F^1=0$. For $i>1$ $\phi_{0i}$ factors through $\phi_{01}$; we obtain $F^i=0$ for $i \geq 1$. Hence from (3), we have $(i-1)!$ $\mathbb{A}^i(X)=0$ and our proof is complete.
\smallskip

\textbf{Remark.} The reader is referred to [P-Sh] for a somewhat different approach for proving equation (5).

\smallskip

We offer a different proof of equation (5) in corollary 4.3.

\smallskip

We have the following corollaries.

\smallskip
\textbf{3.2. Corollary.} Let $(R, m)$ be a ramified regular local ring of dimension $n+1$ essentially of finite type over a discrete valuation ring $(V, pV)$. Then $(i-1)!$ $\mathbb{A}^i(R[1/p])=0$.

Since $V[1/p]$ is a field and $R[1/p]$ is a smooth variety of dimension $n$ over $V[1/p]$, the proof follows readily from the theorem.

\smallskip
\textbf{3.3. Corollary.} If $R^h$ is the henselization of $R$, i.e., $R^h= \varinjlim R_{\lambda}$, where $R_{\lambda}$ is a pointed etale e\'xtension of $R$, then $(i-1)!$ $\mathbb{A}^i(R^h[1/p])=0$.

\smallskip
\textbf{Proof.} Since $R^h[1/p] = \varinjlim R_{\lambda}[1/p]$ and $\mathbb{A}^i (R^h[1/p]) = \varinjlim \mathbb{A}^i (R_{\lambda}[1/p])$, the proof follows from the above corollary via a direct limit argument.

\smallskip
\textbf{3.4 Corollary.} Let $(R, m)$ be as above. Let $\bar R=R/pR$. Then for every prime ideal $P$ of height $i$ such that $p\notin P,$ $(i-1)!$ $[\bar R/\bar P]=0 \in \mathbb{A}^i(\bar R)$.

\smallskip
\textbf{Proof.} We have an exact sequence:
$$\mathbb{A}^{i-1}(\bar R)\rightarrow \mathbb{A}^i(R) \rightarrow \mathbb{A}^{i}(R[1/p]) \rightarrow 0 .......(6).$$
We  also have a group homomorphism

$\theta_i:\mathbb{A}^i(R) \rightarrow \mathbb{A}^i(\bar{R})$ such that if $p \notin P$,  $\theta_i([R/P]) = [\bar{R}/\bar{P}]$ and if $p \in P$, $\theta_i ([R/P]) = 0$ (Gysin map) ......(7). 
\smallskip

Due to (6) and (7), we obtain a group homomorphism $\eta_i: \mathbb{A}^i[R[1/p]] \rightarrow \mathbb{A}^i (\bar{R})$, $\eta_{i} ([R/P]) = [\bar{R} / \bar{P}]$. The desired result now follows from corollary 3.2.

Note here that $\bar{R}$ is not a regular local ring. 

\medskip

\begin{center}{Section 4: Mixed characteristic case and applications}
\end{center}

\smallskip

\textbf{4.1.} In this section we extend Fulton's arguments in chapter 15.3, [Fu] to present Riemann-Roch theorem without denominators over a Dedekind domain in any characteristic. Let us mention briefly the problems that arise in this situation and Fulton's approach to overcome them preceded by works of Samuel, Grothendieck and others in this area.

\smallskip

a) One of the key facts used in studying Chow groups for varieties over a field is the following: in an affine domain $A$ over a field $k$, for any prime ideal $P$, we have dim$A/P+$ ht$P=$ dim$A$. This implies that every maximal ideal of $A$ has the same height $(=$ dim $A)$. Unfortunately, in any characteristic even when $A = V[X_1, ..., X_n]$, $n \geq 1$ and $V$ a dvr, this property fails to hold, i.e., $A$ has maximal ideals of height $n$. In order to avoid such situations Fulton introduced the notion of relative dimension in any characteristic in the following way: let $S$ be an arbitrary regular scheme and let $f: X \rightarrow S$ be a scheme of finite type over $S$. Let $\tilde{V} \subset X$ be a closed integral subscheme; define dim$_S \tilde{V} =$ tr. deg$(R (\tilde{V})/R(T))-$ codim $(T, S)$, where $T=$ closure of $f(\tilde{V}), R(\tilde{V}), R(T)$ are the corresponding function fields. Simply put: if $S=$ Spec$A$, $X=$ Spec$B$, $\tilde{V}=$ Spec$(B/P)$, and $P \cap S = q$, then dim$_S \tilde{V}=$tr$_{k(q)}k(P)-$ ht$q$. For example, if $A=V[X_1, ..., X_n]$, $V$ a dvr, $m$ is a maximal ideal of height $n+1$, then dim$_V A/m =-1$ and if $m$ is a maximal ideal of height $n$ then dim$_V A/m =0$. This immediately shows that 

i) for $U$ non-empty open in $\tilde{V}$, dim$_SU=$ dim$_S \tilde{V}$

ii) $\tilde{V}$ as above; dim$_SX=$ dim$_S \tilde{V} + $ codim $(\tilde{V}, X)$ and

iii) if $f:\tilde{V} \rightarrow W$ is a dominant morphism of varieties over $S$ then dim$_S \tilde{V}=$ dim$_SW+$tr. deg$(R (\tilde{V})/R (W))$. 

\smallskip

These three properties pave the way for defining t-cycles as $\Sigma n_i V_i, V_i$- a subvarity of $X$, dim$_S V_i=$t and rational equivalence etc. We denote the corresponding group modulo rational equivalence by $\mathbb{A}_t (X/S)$ or $\mathbb{A}_t (X)$ when there is no scope for confusion. Gr$(K_0(X))$ is also defined by taking into consideration the filtration on $K_0(X)$ obtained by using the notion of relative dimension. We would urge the reader to go through p.394-395 in [Fu] for a brief description of related facts about Chern classes, Chern characters, e.g., for any vector bundle $E$ on $X$ we have, $c_i(E):\mathbb{A}_t(X/S) \rightarrow \mathbb{A}_{t-i}(X/S)$, their corresponding properties and major theorems involving Chern classes and Chern characters in intersection theory.
\medskip

b) However, even with all these properties, exterior product of Chow groups does not work in general, i.e., one can not define $\mathbb{A}_i (X/S) \otimes \mathbb{A}_j (Y/S) \rightarrow \mathbb{A}_{i+j} (X \underset{S}{\times} Y/S)$ (subvarities of $X$ may not be flat over $S$).
\smallskip

In order to remove this obstacle Fulton restricted $S$ to Spec$A$ where $A$ is a Dedekind domain. This means any variety over $S$ is either flat or maps to a closed point in $S$. In this situation, given $V \subset X, W \subset Y$ subvarities, Fulton defines a product cycle $[V] \underset{S}{\times} [W]$ on $X \underset{S}{\times} Y$ as follows:
\smallskip

$[V] \underset{S}{\times} [W] = \{[V \underset{S}{\times} W]$, if $V$ or $W$ is flat over $S$; $0$ otherwise $\}$

and proves that this product respects rational equivalence and defines an exterior product 
\smallskip

$\mathbb{A}_i (X/S) \otimes \mathbb{A}_j (Y/S) \rightarrow \mathbb{A}_{i+j} (X \underset{S}{\times} Y/S)$, (prop. 20.2, [Fu]). This construction satisfies usual properties (e.g., commutativity, associativity etc.) for exterior products; for schemes smooth over $S$ of relative dimension $n$, intersection multiplicity, ring structure on $\mathbb{A}^*(X/S)=\oplus \mathbb{A}^i(X/S), \mathbb{A}^i (X/S) = \mathbb{A}_{n-i}(X/S)$ (Ch. 8, [Fu]) make sense without any change in corresponding properties (e.g., product via diagonal imbedding, graph construction via $f: X \rightarrow Y$, a morphism of smooth scheme over $S$ etc.). In short, whatever operations we perform on a smooth variety over a field can be carried over to a variety smooth over $S$ of relative dimension $n$ where $S=$ Spec $A$, $A$- a Dedekind domain. We refer to (20.2) in [Fu] for details. 
\smallskip

In particular we would like to mention the following observation which Fulton proved using bivariant property of Riemann-Roch: 
\smallskip

\textbf{Observation.} If $X$ is smooth of relative dimension $n$ over $S$, $S$ as above, then $\mathbb{A}^i(X) \simeq \mathbb{A}_{n-i}(X/S)=\mathbb{A}^i(X/S)$.

This will be used in our work. 
\medskip

{\bf 4.2.} Now we are ready to state and sketch a proof of Riemann-Roch theorem without denominators in any characteristic over a Dedekind domain extending the steps described in (15.3, [Fu]). For our proof sometimes we highlight part of Fulton's proof leaving the details to the reader and sometimes we elaborate part of his approach that he sketched very briefly. In this respect we only extend the arguments for crucial steps in Fulton's proof. First we mention the following observation.

\smallskip

\textbf{Observation 1.} (Lemma 15.3, [Fu]) For fixed positive integers $d, e$ there is a unique power series with integer coefficients $g(T_1, ..., T_d, U_1, ..., \\ U_e)$ such that for all vector bundles $Q, E$ of ranks $d, e$ respectively on any variety $X$ smooth of relative dimension $n$ over $S$, $S=$ spec$A$, $A$ - a Dedekind domain. we have 

$$c(\land^{\bullet} Q^{\lor} \otimes E) -1 = c_d(Q).g(Q, E)$$
where $g(Q, E)$ denotes $g(c_1 (Q), ..., c_d (Q), c_1 (E), ..., c_e (E))$.
\smallskip

\textbf{Sketch of a proof.} Due to the validity of theorem 3.2 (e) (Whitney Sum), Remark 3.2.3 (c) (exterior power) and example 3.2.5 (Lemma 18, [B-Se]) in [Fu] for relative set up in any characteristic over Dedekind domains, proof of this observation is obtained by the same arguments as in the proof in [Fu] for Lemma 15.3.
\smallskip

\textbf{Theorem.} (Riemann-Roch theorem without denominators in any characteristic, [J],  Theorem 15.3, [Fu]) Let $f: X \rightarrow Y$ be a closed imbedding of smooth varieties over $S=$ Spec$A$, $A$ - a Dedekind domain in any characteristic, of codimension $d$ with normal bundle $N$. Let $E$ be a vector bundle of rank $E$ on $X$. Then,
\smallskip

\ \ $c(f_* [E])=1+f_* (g(N, E))$, where $g$ is defined as in observation 1.
\smallskip

\textbf{Sketch of a proof.} As we mentioned earlier, Fulton's proof for theorem 15.3 in equicharacteristic works here without any scope for confusion. We want to describe here a special case of the above theorem (model for closed imbedding) in our situation to point out the utility of observation 1.

\smallskip
Let $Y=\mathbb{P} (N \oplus 1)$, where $N$ is a vector bundle of rank $d$ on $X$; $f=$ the composite of $X \hookrightarrow N \hookrightarrow \mathbb{P} (N \oplus 1)$ where $X \hookrightarrow$ zero-section of $N$, $N \hookrightarrow \mathbb{P} (N \oplus 1)$ is the canonical open imbedding. Let $\pi:\mathbb{P} (N \oplus 1) \rightarrow X$ be the projection map. We have an exact sequence
$$0 \rightarrow O_Y(-1) \rightarrow \pi^*(N \oplus 1) \xrightarrow{\eta} Q \rightarrow 0    ............(7)$$
where $Q$ is the universal quotient  bundle. Let $s=$ the section of $Q$ defined by $\eta(\pi^*(1))$. Then $Z(s)(=$ the zero-set of $s)$ is exactly $X$. Hence $c_d(Q).[Y]=f_*[X]$, for any $\alpha \in \mathbb{A}^*(Y/S), f_*(f^* \alpha)=\alpha.f^*[X]=c_d(Q).\alpha .................(8)$. 

Due to the fact that $X=Z(s)$, it follows that $f_* O_X$ has a resolution determined by the Koszul complex defined by $s:$

$$0\rightarrow \overset{d}{\land} Q^{\lor} \rightarrow ... \rightarrow \overset{2}{\land} Q^{\lor} \rightarrow Q^{\lor} \xrightarrow{s^{\lor}} O_Y\rightarrow f_* O_X \rightarrow 0 ....(9)$$
This implies that $f_* E$ has the following resolution.
$$0\rightarrow \overset{d}{\land} Q^{\lor} \otimes \pi^* E \rightarrow ... \rightarrow \overset{2}{\land} Q^{\lor} \otimes \pi^* E \rightarrow Q^{\lor} \otimes \pi^* E \xrightarrow{s^{\lor}} \pi^* E \rightarrow f_* E \rightarrow 0.$$
Hence, by observation 1, c$(f^* [E]) = 1+$c$_d(Q).g(Q, \pi^* E)$

Since $f^*(Q)=N, f^* \pi^*(E)=E$, due to (8) it follows that c$ (f_* [E])= 1+f_*(g(N, E))$
\smallskip

In his proof of Grothendieck's Riemann-Roch theorem (Theorem 15.2, [Fu]) Fulton uses deformation to the normal bundle to deform $f: X \rightarrow Y$, a closed imbedding, $X, Y$ as above, into the imbedding $\tilde{f}: X \rightarrow \mathbb{P}(N \oplus 1), N$ being the normal bundle to $X$ in $Y$. Fulton pointed out that the same technique works for his proof of Riemann-Roch theorem without denominators (theorem 15.3 in [Fu]) by changing Chern characters to Chern classes in the above proof. This technique works as well for any characteristic case without any scope for confusion (replace $X \times \mathbb{P}^1$ by $X \underset{S}{\times}\mathbb{P}_S^1$).
\medskip

\textbf{4.3 Corollary.} (eqn.(16), p.151, [Gr1]; [Gr2]; 15.3.6, [Fu]; [P-Sh]) With hypothesis and notations as above we have, 

c$_j(f_*[O_X])=0$ for $0<j<d$, c$_d(f_*[O_X])=(-1)^{d-1}(d-1)!$ $[X]$ in $\mathbb{A}^d (Y/S)$.
\smallskip

\textbf{Proof.} We could derive the proof using the above theorem (Example 15.3.1 in [Fu]). We here provide a different proof for all characteristics over a Dedekind domain.

we have an exact sequence.
$$\mathbb{A}_k(X/S)\rightarrow \mathbb{A}_k (Y/S)\rightarrow \mathbb{A}_k(Y-X/S) \rightarrow 0 ........(10)$$

If dim$_SX<k$ then $\mathbb{A}_k (Y/S) \simeq \mathbb{A}_k(Y-X/S)$.

It follows from (9) that 

c$(f_*[O_X])=$c$(\Sigma(-1)^i \overset{i}{\land} Q^{\lor})= \Pi$ c $(\overset{i}{\land}Q^{\lor})^{(-1)^i}$. Then, for any $j>0$, $c_j(f_*[O_X])=0$ on $Y-X$. Hence, it follows from (10) that for $j<d$, 

\smallskip
$c_j (f_*[O_X])=0$ and $c_d (f_* [O_X])=$ a multiple of 
$[X] \in \mathbb{A}^d(Y/S) ....(11)$. 

\smallskip
Let us recall that for any vector bundle $E$ of rank $r$, if $\alpha_1, ..., \alpha_r$ denotes the Chern roots of $E$ and $p_k=\alpha_1^k+ ... +\alpha_r^k$, then $ch(E)=\Sigma1/k!p_k$.

Let $c_1, c_2, ..., c_k, ...$ denote the Chern classes of $E$.

Macdonald has shown that (p.30, [M]).

%Let us recall that ch$(f_*(O_X))=\Sigma (1/k!)p_k$ where 

\medskip

%\begin{equation}
$p_k=det 
\begin{pmatrix} c_1 & 1 & 0 & ... &  0 &  0 \\
                2c_2 & c_1 & 1 & ... & ... & 0 \\
                3c_3 & 2c_2 & c_1 & 1 & ... \\  
                ... \\
                ... & ... & ... & ... & ... &  1 \\
                kc_k & (k-1)c_{k-1} & ... & ... & ... & c_1 \\
\end{pmatrix}                               ........(12)
%\end{equation}
$

\medskip

It follows from above that $ch(E)=r+c_1+1/2(c_1^2-2c_2)+...$ \ \ (13)

(See also A.IV, p.70-71, [Bou])

If $c_i=0$ for $0<i<d$, then $p_d=(-1)^{d-1}dc_d$ and $ch(E)=r+(-1)^{d-1}/((d-1)!)c_d(E)+...$.        \ \  \ (14)

When $X$ is smooth, due to the additive property of the rank function: $K_0(X) \rightarrow Z$, observations in (12), (13) and (14) extend to every element of $K_0(X)$.

\smallskip
Thus, by (11), (12), (13), (14), we obtain  

\medskip

%Since c$_1 = ... =$ c$_{d-1}=0$, we have $p_1= ... =p_{d-1}=0$ and $p_d = (-1)^{d-1} d$ c$_d$. Hence

ch$(f_*(O_X))=(-1)^{d-1}d c_d/d!$ $+p_{d+1}/(d+1)!$ $+...$

\ \ \ \ \ \ \ \ \ \ \ \ \ \ \ \ \  $=(-1)^{d-1} c_d/(d-1)!$ $+...$

By Riemann-Roch theorem, we have 

ch$(f_*(O_X)) .$ td$(T_Y)=f_*($ch $O_X .$ td$(T_X))$

Equating the initial terms, from above, we obtain $(-1)^{d-1}/(d-1)!$ $c_d (f_* [O_X]) = [X]$, i.e., $c_d(f_*[O_X])=(-1)^{d-1}(d-1)!$ $[X]$ in $\mathbb{A}^d(Y/S)_{\mathbb{Q}}$.

Hence the assertion follows from (11).
\medskip

Our final theorem in this section is a  generalised version of Theorem 3.1

\medskip

\textbf{4.4 Theorem.} Let $X$ be a variety of dimension $n+1$ smooth over spec$A$, $A$- a Dedekind domain. Assume that $K_0 (X)= \mathbb{Z}$;  $K_0 (X)= K_0 (\mathscr{M}_0)$ where $\mathscr{M}_0$ is the category of coherent sheaves on $X$. Then $(i-1)!$ $\mathbb{A}^i(X)=0$, for $0 <i <n+1$.
\medskip

\textbf{Proof.} Due to corollary 4.3 our proof is obtained by following exactly the same arguments as were used in the proof of Theorem 3.1. (replacing $\mathbb{A}^i(X)$ by $\mathbb{A}^i(X/S)$).
\medskip

We have the following corollaries:
\smallskip

\textbf{4.5 Corollary.} Let $(R, m)$ be a ramified regular local ring of dimension $(n+1)$ essentially finite type over a pseudo-geometric dvr $(V, pV)$ such that $V/pV$ is perfect or $V/pV = R/m$, i.e., $R=(S[X]/ \\ (f(X)))_{\tilde{m}}$ by theorem 2.1. Then the unramified locus of Spec$R$ over Spec$S$, i.e., Spec$R_{f'}$ has the property that, for $0<i \leq n, (i-1)!$ $\mathbb{A}^i(R_{f'})=0$.
\smallskip

Recall from theorem 2.1 that $S$ is a localization of a polynomial ring over a dvr $V$ at a maximal ideal. Hence, $V \rightarrow S$ is smooth and $S \rightarrow R_{f'}$ is \'etale; thus, $V \rightarrow R_{f'}$ is smooth. Now the assertion follows from the above theorem.

\smallskip

\textbf{4.6 Corollary.} Let $(R, m)$ be as above and let $(R^h, m^h)$ be its henselization. Let $U^h=$ quasi-unramified locus of spec$R^h$ over spec$S$. Then $(i-1)!$ $\mathbb{A}^i(U^h) =0$, for $0<i \leq n$.
\smallskip

Since $R^h=\varinjlim R_{\lambda}$ where $R_{\lambda}$ is a pointed \'etale extension of $R$, $U^h=\varinjlim$ of $U_{\lambda}, U_{\lambda}=$ unramified locus of spec$R_{\lambda}$ over spec$S$. Our result now follows from the above corollary via a direct limit argument (cor.3.3).
\smallskip

\textbf{4.7 Corollary.} Let $(R, m)$ be as in corollary above. Let $P$ be a prime ideal of $R$ such that ht$\bar{P=i}$ and $f' \notin P$. Let $\bar R=R/f'R$ and $\bar R/\bar P = R/(P+f' R)$. Then $(i-i)!$ $[\bar R/\bar P]=0$ in $\mathbb{A}^i_{\bar{R}}$.
\smallskip

The proof is identical with that of corollary 3.4 replacing $p$ by $f'$.

\medskip

\begin{center} Section 5: Relation between $\mathbb{A}^i (R^h)$ and $\mathbb{A}^i \hat{R}$ \end{center}

\textbf{5.1 A very brief preview of Artin approximation and related results.}  

\smallskip
A Noetherian local ring $(A, m)$ is called an approximation ring if the following condition holds: if any system of polynomial equations $\{F_i(Y_1, ..., Y_n) = 0, 1 \leq i \leq t \}$ over $A$ has a solution $Y_1 = b_1, ..., Y_n = b_n$ in the $m$-adic completion $\hat{A}$ of $A$, then, given any positive integer $N$, this system must have a solution $Y_1 = a_1, ..., Y_n = a_n$ in $A$ such that $a_i \equiv b_i$ mod $m^N \hat{A}$ for $1 \leq i \leq n$.

\smallskip
Actually this is equivalent to the statement that whenever a systen of polynomial equations over $A$ has a solution in $\hat{A}$, the system has a solution in $A$. For most applications we need the stronger version as stated above. It is easy to see that an approximation ring is Henselian. Artin's notable contribution is the following:

\smallskip
\textbf{Theorem: (Artin approximation)} let $S$ be a finitely generated algebra over $K$, where $K$ is a field or an excellent dvr. Let $R = S_P, P$ is a prime ideal in $S$. Then $R^h$ (the Henselization of $R$) is an approximation ring. Also $\mathbb{C} \{X_1, ..., X_n\}/I$ is an approximation ring.

This theorem has a tremendous impact on research in commutative algebra and algebraic geometry. Several homological questions / conjectures on complete local rings have been reduced to the same on local rings essentially of finite type over a field or an excellent dvr via this process of approximation. For results relevant to this paper we refer the reader to [P-S 1] and [D3]. In particular we want to mention the following results.

\smallskip

Let $V$ be a field or an excellent dvr and $A$ be a local ring of essentially finite type over $V$. Let $\hat{A}$ denote the completion of $A$ with respect to maximal ideal $m$ of $A$ and $A^h$ denote the $m$-adic henselization of $A$ in $\hat{A}$. Then $A^h = \underset{n \in I}{\varinjlim}A_n$, where each $A_n$ s a pointed \'etale extension of $A$. 

\smallskip

{\bf Result 1.} (Proposition 3.3 in [D3]) Let $0 \rightarrow \hat{N} \rightarrow \hat{M} \rightarrow \hat{T} \rightarrow 0$ be an exact sequence  of finitely generated $\hat{A}$-modules and let $\hat{F}_\bullet$, $\hat{L}_\bullet$ $\hat{G}_\bullet$ be free resolution of $\hat{N}$, $\hat{M}$, $\hat{T}$ respectively where $\hat{F}_\bullet$, $\hat{G}_\bullet$ are minimal and $\hat{L}_\bullet = \hat{F}_\bullet \oplus \hat{G}_\bullet$ is constructed from $\hat{F}_\bullet, \hat{G}_\bullet$ in the usual manner i.e. if $\hat{\phi}_\bullet, \psi_\bullet$ represents the boundary maps of $\hat{F}_\bullet, \hat{G}_\bullet$ respectively; and $\hat{h}_\bullet$ represents the map from $(\hat{G}_\bullet)_1 \rightarrow F_\bullet$ with $\hat{h}_n \hat{\psi}_{n+1} = 0$, then $\hat{\theta}_\bullet =
  \begin{pmatrix}
    \hat{\theta}_\bullet & \hat{h}_\bullet\\
    0 & \hat{\psi}_\bullet\\
  \end{pmatrix}$
represents the boundary maps of $L_\bullet$. Then for every positive integer $t$, there exists an $n \in I$ and an exact sequence of finitely generated $A_n$-modules
$$0 \rightarrow N_n \rightarrow M_n \rightarrow T_n \rightarrow 0$$
such that

\smallskip
(i) $N_n \otimes A_n/m^t_n = \hat{N} \otimes \hat{A}/ \hat{m}^t$, $M_n \otimes A_n/m^t_n = \hat{M} \otimes \hat{A}/ \hat{m}^t$, $T_n \otimes A_n/ m^t_n = \hat{T} \otimes \hat{A}/ \hat{m}^t$.

\smallskip
(ii) There exists free resolutions $F_\bullet$, $L_\bullet$, $G_\bullet$ of $N_n$, $M_n$ and $T_n$ respectively such that

\smallskip
    (a) $\hat{F}_\bullet$, $\hat{G}_\bullet$ are minimal and $L_\bullet = F_\bullet \oplus G_\bullet$ is constructed in the usual manner,

\smallskip
    (b) $F_\bullet \oplus A_n/m^t_n = \hat{F}_\bullet \otimes \hat{A}/ \hat{m}^t$, $L_\bullet \otimes A_n/ m^t_n = \hat{L}_bullet \otimes \hat{A}/ \hat{m}^t$ and $G_\bullet \otimes A_n/ m^t_n = \hat{G}_\bullet \otimes \hat{A}/ \hat{m}^t$.    

\bigskip
\textbf{Result 2.} (Proposition 3.5 [D3]). Let $\hat{M}$ be a finitely generated $\hat{A}$-module with dim$\hat{M} = d$. Then there exists a finitely generated module $M_n$ over $A_n$ for some $n \in I$ such that dim $M_n = d$ and $M_n \otimes \hat{A}/ \hat{m}^t = \hat{M}/m^t \hat{M}$.

\smallskip

We apply these results to establish relations between Chow groups of the henselization $(R^h, m^h)$ of a ramified regular local ring $(R, m)$ essentially of finite type over an excellent dvr and its $m$-adic completion $(\hat{R}, \hat{m})$. As stated in the introduction, the definition of $W_i (R)$ in [C-F] converges with the definition of $\mathbb{A}^i (R)$ over any regular local ring $R$. We do not consider relative Chow groups here. Let us recall that for any Noetherian ring $A$, $\mathscr{M}_i(A)$ denotes the Serre-subcategory of finitely generated $A$-modules of codemnesion $\geq i$. We drop $A$ from this notation when there is no scope for confusion. First we prove the following lemma. 

\smallskip

{\bf 5.2 Lemma.} $\mathbb{A}^i(R)=0$ for $1 \leq i \leq n \Longleftrightarrow K_0(\mathscr{M}_i) \xrightarrow {\phi_{i-1, i}} K_0(\mathscr{M}_{i-1})$ are $0$-maps for all $i$, $1 \leq i \leq n$. This part is valid for any regular local ring.

\smallskip

{\bf Proof.} We have the following commutative diagram of exact sequences:
\medskip

%\begin{equation}\label{CD2}
\begin{tikzcd}[cramped, sep=small]
& K_0(\mathscr{M}_{i+1}) \arrow[equal]{r} \arrow[d,"\phi_{i, i+1}"] & K_0(\mathscr{M}_{i+1}) \arrow[d, "\Phi_{i-1, i+1}"] \\
&  K_0(\mathscr{M}_i) \arrow[r, "\phi_{i-1, i}"] \arrow[d,"\eta_i"] & K_0(\mathscr{M}_{i-1}) \arrow[r, "\eta_{i-1}"] \arrow[d, "\tilde{\eta}_{i-1}"] & K_0(\mathscr{M}_{i-1}/ \mathscr{M}_i) \arrow[r] \arrow[equal]{d} &  0  & (*)\\
& K_0(\mathscr{M}_i/ \mathscr{M}_{i+1}) \arrow[r, "\bar{\phi}_{i-1, i}"] \arrow[d] & K_0(\mathscr{M}_{i-1}/ \mathscr{M}_{i+1}) \arrow[r, "\bar{\eta}_{i-1}"] \arrow[d] & K_0(\mathscr{M}_{i-1}/ \mathscr{M}_i) \arrow[r] & 0 \\
& 0 & 0 
\end{tikzcd}
%\end{equation}
\smallskip

where the vertical columns are also exact.
\smallskip

We have $\mathbb{A}^i(R)$ = im $\bar{\phi}_{i-1, i}$ (p.230, [C-F]).
From (*) it follows that $\phi_{i-1, i}=0 \Rightarrow \mathbb{A}^i(R) = 0$. Conversely, suppose $\mathbb{A}^i(R)=0$ for $0 < i \leq n$. We induct on $i$. Since $R$ is regular local $K_0(\mathscr{M}_n) \xrightarrow{\phi_{n-1, n}} K_0(\mathscr{M}_{n-1})$ is the 0-map.
\smallskip

Hence from (*) we obtain 

\newcommand*{\isoarrow}[1]{\arrow[#1,"\rotatebox{90}{\(\sim\)}"
]}

\[
\begin{tikzcd}
    & K_0(\mathscr{M}_{n-1})\arrow{r}{\phi_{n-2, n-1}}\isoarrow{d} & K_0({\mathscr{M}_{n-2}})\isoarrow{d} \\
    & K_0(\mathscr{M}_{n-1}/\mathscr{M}_n)\arrow{r}{\bar{\phi}_{n-2, n-1}} & K_0(\mathscr{M}_{n-2}/\mathscr{M}_n) \\
\end{tikzcd}
\]

Since $\mathbb{A}^{n-1}(R)=0 \Rightarrow \bar{\phi}_{n-2, n-1}=0$, we have $\phi_{n-2, n-1}=0$ (from the above diagram).

Suppose we have proved that $K_0(\mathscr{M}_j) \xrightarrow{\phi_{j-1, j}} K_0(\mathscr{M}_{j-1})$ are $0$-maps for $j>i$. Hence from (*) we obtain vertical isomorphism

\[
\begin{tikzcd}
 & K_0(\mathscr{M}_i) \arrow{r} {\phi_{i-1, i}} \isoarrow{d} & K_0({\mathscr{M}_{i-1}}) \isoarrow{d} \\ 
 & K_0(\mathscr{M}_i/\mathscr{M}_{i+1}) \arrow{r}{\bar{\phi}_{i-1, i}}  & K_0(\mathscr{M}_{i-1}/\mathscr{M}_{i+1}) \\ 
\end{tikzcd}
\]

Since $\mathbb{A}^i(R) = 0 \Rightarrow \bar{\phi}_{i-1, i} = 0$, we have $\phi_{i-1, i}=0$.

\smallskip

\textbf{Remark.} The above proof shows that the validity of the $K_0$-case of Gersten's conjecture [Q] is equivalent to the validity of the Chow group problem over any regular local ring.

\smallskip

Next we prove the following proposition.

{\bf 5.3 Proposition.} For $1 \leq i \leq n, K_0(\mathscr{M}_i(R^h)) \rightarrow K_0(\mathscr{M}_i(\hat{R}))$ is injective.

\smallskip

{\bf Proof.} Let $M, N \in \mathscr{M}_i(R^h)$ be such that $[M \otimes_{R^h} \hat{R}] = [N \otimes_{R^h} \hat{R}]$ in $K_0(\mathscr{M}_i (\hat{R}))$. Then there exist short exact sequences

$0 \rightarrow {T'} \rightarrow T \rightarrow T'' \rightarrow 0$ and $0 \rightarrow W' \rightarrow W \rightarrow W'' \rightarrow 0$, $T$, $W$ etc. $\in (\mathscr{M}_i(\hat{R}))$ such that

$(\mathscr{M} \otimes_{R^h} \hat{R}) \oplus T' \oplus T'' \oplus W \simeq (N \otimes_{R^h} \hat{R}) \oplus T \oplus W' \oplus W''$  (p.74, [Sw])

Then, by Result 1 and Result 2 mentioned in 5.1 and related results in [D3], there exist short exact sequences

$0 \rightarrow U' \rightarrow U \rightarrow U'' \rightarrow 0$ and $0 \rightarrow L' \rightarrow L \rightarrow L'' \rightarrow 0 \in \mathscr{M}_i(R^h)$ such that

$M \oplus U' \oplus U'' \oplus L \simeq N \oplus U \oplus L' \oplus L''$

Hence $[M] = [N] \in K_0(\mathscr{M}_i(R^h))$.

\smallskip

Now we are ready to state our theorem and a corollary.

{\bf 5.4 Theorem.} If $\mathbb{A}^j(\hat{R})=0$ for $i \leq j \leq n$, then $\mathbb{A}^j(R^h)=0$ for $i \leq j \leq n$.

\smallskip

{\bf Proof.} Suppose $\mathbb{A}^j(\hat{R})=0$ for $i \leq j \leq n$. Then, by the proof of Lemma 5.2, we have

$K_0(\mathscr{M}_j(\hat{R})) \xrightarrow{\hat{\phi}_{j-1, j}} K_0(\mathscr{M}_{j-1}(\hat{R}))$ are $0$-maps for $j\geq i$. We have a commutative diagram
\smallskip

\begin{tikzcd}
    & K_0(\mathscr{M}_j(R^h))\arrow{r}{\phi^h_{j-1, j}}\arrow{d}{\Theta_j} & K_0({\mathscr{M}_{j-1}}(R^h))\arrow{d}{\Theta_{j-1}} \\
    & K_0(\mathscr{M}_j(\hat{R}))\arrow{r}{\hat{\phi}_{j-1, j}} & K_0(\mathscr{M}_{j-1}(\hat{R}) \\
\end{tikzcd}

where $\Theta_j$s are induced by the flat map $R^h \rightarrow \hat{R}$. By the above proposition, $\Theta_j$s are injective. Since $\hat{\phi}_{j-1, j}=0$, $\phi^h_{j-i, j}$ are $0$-maps for $j \geq i$. Hence, by Lemma 5.2, $\mathbb{A}^j(R^h)=0$ for $j \geq i$.

\smallskip

{\bf 5.5 Corollary} $\mathbb{A}^{n-1}(R^h)=0$.

\smallskip

{\bf Proof.} Since $R^h$ and $\hat{R}$ are both regular local, $0=\mathbb{A}^n(R^h)=\mathbb{A}^n(\hat{R})$. By theorem 1.1 in [D1], we have $\mathbb{A}^{n-1}(\hat{R})=0$. Hence, by the above theorem, $\mathbb{A}^{n-1}(R^h)=0$.

\smallskip

\end{document}